\documentclass[12pt,a4paper]{article}
\usepackage[cp1251]{inputenc}
\usepackage[T2A]{fontenc}
\usepackage[english]{babel}

\usepackage{amsmath,amsfonts,amssymb,mathrsfs,amsopn,indentfirst,amscd}

\usepackage{graphicx}
\usepackage{amsthm}
\usepackage{hyperref}

\usepackage{amssymb}
\usepackage{amsmath}
\newcommand{\ver}{{\rm ver}}
\newcommand{\vo}{{\rm vol}}

\newtheorem*{corollary*}{Corollary}

\begin{document}
\title{Properties of 0/1-Matrices of Order  $n$ Having Maximum Determinant
}
\author{Mikhail Nevskii\footnote{Department of Mathematics,
              P.G.~Demidov Yaroslavl State University, Sovetskaya str., 14, Yaroslavl, 150003, Russia 
              orcid.org/0000-0002-6392-7618 
              mnevsk55@yandex.ru}        \and
        Alexey Ukhalov \footnote{ Department of Mathematics, P.G.~Demidov Yaroslavl State University, Sovetskaya str., 14, Yaroslavl, 150003, Russia 
              orcid.org/0000-0001-6551-5118 
              alex-uhalov@yandex.ru}
              }
\date{July 4, 2018}
\maketitle

\begin{abstract}
We give some necessary conditions for maximality of 0/1-determinant.
Let ${\bf M}$ be a  nondegenerate $0/1$-matrix of order $n$.
Denote by $\bf A$  the matrix of order $n+1$ which appears
from ${\bf M}$
after adding the $(n+1)$th row $(0,0,\ldots,0,1)$ and the $(n+1)$th column
consisting of $1$'s. Suppose ${\bf A}^{-1}=(l_{ij}),$  then
 for all $i=1,\ldots,n$
we have
$\sum_{j=1}^{n+1} |l_{ij}|\geq 2.$ 
Moreover, if $|\det({\bf M})|$ is equal to the maximum value of a  $0/1$-determinant 
of order $n$, then $\sum_{j=1}^{n+1} |l_{ij}|= 2$ for all $i=1,\ldots,n$.
\medskip 

\end{abstract}

\section{Introduction}
\label{nev_ukh_intro}

In this paper we present some properties of 0/1-matrices. 
Our approach to achieve these results lies in the field between linear algebra and geometry of convex bodies. At the same time, we demonstrate new examples how geometric results can be expressed 
by means of linear algebra. This is why we start with necessary geometric definitions.

Assume $n\in{\mathbb N}$, 
$Q_n:=[0,1]^n$.
Let $S$ be a nondegenerate simplex $S\subset{\mathbb R}^n$.  
Denote by $\sigma S$ the homothetic copy of $S$ with center of homothety 
at the center of gravity of $S$ and ratio $\sigma$. 
By symbol $d_i(S)$ we mean  {\it the $i$th axial diameter of  $S$}, i.\,e. 
the length of a longest segment in $S$ parallel to the $i$th coordinate axis.
The notion of axial diameter of a convex body was introduced by Paul Scott  
(see \cite{nu_bib_10} and 
\cite{nu_bib_11}
).
Let us define $\alpha(S)$  as the minimal $\sigma>0$ such that
$Q_n$ is a~subset of a translate of $\sigma S$. Consider the value
$\xi(S):=\min \{\sigma\geq 1: Q_n\subset \sigma S\}.$
Obviously, for any simplex $S$, we have $\alpha(S)\leq \xi(S)$. By definition,
put
$$
\xi_n:=\min \{ \xi(S): \, S \mbox{ --- $n$-dimensional simplex,} \, 
S\subset Q_n, \, \vo(S)\ne 0\}.
$$

Remark that all the values $d_i(S)$, $\alpha(S)$, and $\xi(S)$ can be calculated with
the use of vertices of $S$. 
Denote these vertices 
by $x^{(j)}=\left(x_1^{(j)},\ldots,x_n^{(j)}\right),$ 
$1\leq j\leq n+1,$ and define the following {\it node matrix of the  simplex $S$:} 
$$
{\bf A} :=
\left( \begin{array}{cccc}
x_1^{(1)}&\ldots&x_n^{(1)}&1\\
x_1^{(2)}&\ldots&x_n^{(2)}&1\\
\vdots&\vdots&\vdots&\vdots\\
x_1^{(n+1)}&\ldots&x_n^{(n+1)}&1\\
\end{array}
\right).
$$
Let $\Delta:=\det({\bf A})$, then $\vo(S)=\frac{|\Delta|}{n!}$.
Denote by $\Delta_j(x)$ the determinant  obtained from 
 $\Delta$ by changing the  $j$th row with  the row
$(x_1,\ldots,x_n, 1).$ 
Linear polynomials
$\lambda_j(x):=
\frac{\Delta_j(x)}{\Delta}$ 
satisfy the property
$\lambda_j\left(x^{(k)}\right)$ $=$ 
$\delta_j^k$, where $\delta_j^k$ is the  Kronecker delta-symbol. 
Coefficients of $\lambda_j$ form the $j$th column of ${\bf A}^{-1}$.
Assume
${\bf A}^{-1}$ $=(l_{ij})$, i.\,e.,~$\lambda_j(x)=
l_{1j}x_1+\ldots+
l_{nj}x_n+l_{n+1,j}.$

Since each linear polynomial  $p$ 
can be represented in the form
$$p(x)=\sum_{j=1}^{n+1} p\left(x^{(j)}\right)\lambda_j(x),$$
we call  $\lambda_j$ 
{\it basic Lagrange polynomials related to $S.$}
By taking  $p(x)=x_1, \ldots, x_n, 1$, we get the equalities
\begin{equation} \label{nev_uhl_lagr_sums}
\sum_{j=1}^{n+1} \lambda_j(x) x^{(j)}=x, \quad
\sum_{j=1}^{n+1} \lambda_j(x)=1.
\end{equation}
Consequently, $\lambda_j(x)$  
are the barycentric coordinates of a point $x\in{\mathbb R}^n$ with respect to 
$S$. The simplex $S$ can be determined by the system of relations
$\lambda_j(x)\geq 0.$ 

The $i$th axial diameter of $S$ satisfies the equality
\begin{equation}\label{nev_uhl_d_i_formula}
\frac{1}{d_i(S)}=\frac{1}{2}\sum_{j=1}^{n+1} \left|l_{ij}\right|. 
\end{equation}
There exist only one line segment in $S$ with the length $d_i(S)$ parallel to
the $i$th coordinate axis. Each $(n-1)$-facet of $S$ 
contains at least one of the endpoints of this segment.
The center $c^{(i)}$ of the segment can be expressed through the numbers
$l_{ij}$ and vertices $x^{(j)}$:
\begin{equation}\label{nev_uhl_max_center}
c^{(i)}= \sum_{j=1}^{n+1} m_{ij} 
x^{(j)}, \quad 
m_{ij}:=
\frac{\left|l_{ij}\right|}
{\sum\limits_{k=1}^{n+1}\left|l_{ik}\right|}.
\end{equation}                                        
If $Q_n\not\subset S$, then we have
\begin{equation}\label{nev_uhl_xi_s_cub_formula}
\xi(S)=(n+1)\max_{1\leq j\leq n+1}
\max_{x\in \ver(Q_n)}(-\lambda_j(x))+1,
\end{equation}
where $\ver(Q_n)$ denotes the set of vertices of $Q_n.$
The equalities
(\ref{nev_uhl_d_i_formula})--(\ref{nev_uhl_xi_s_cub_formula})
were obtained in
\cite{nu_bib_3}.
It was proved in \cite{nu_bib_4}
that
\begin{equation}\label{nev_uhl_alpha_d_i_formula}
\alpha(S)
=\sum_{i=1}^n\frac{1}{d_i(S)}.
\end{equation}
Being combined with (\ref{nev_uhl_d_i_formula}), formula 
(\ref{nev_uhl_alpha_d_i_formula}) implies   the following simple connection between
$\alpha(S)$ and elements of ${\bf A}^{-1}$:
\begin{equation}\label{nev_uhl_alpha_lij_formula}
\alpha(S)=\frac{1}{2}\sum_{i=1}^n\sum_{j=1}^{n+1} |l_{ij}|. 
\end{equation}

If $S\subset Q_n$, then $d_i(S)\leq 1$. Applying (\ref{nev_uhl_alpha_d_i_formula}),
we get under this assumption
$\xi(S)\geq \alpha(S)\geq n.$
Consequently,
$\xi_n\geq n.$  
Various properties and estimates 
of the introduced characteristics are collected in 
\cite{nu_bib_5}.
Later on
some of these estimates
were improved or supplemented in 
\cite{nu_bib_6},
\cite{nu_bib_7},
\cite{nu_bib_8},
\cite{nu_bib_9}.  
Let us mention here some results of the authors concerning $\xi_n.$

It is clear that $\xi_1=1.$ The exact value of $\xi_2$ is equal to
$\frac{3\sqrt{5}}{5}+1=2.34\ldots$
If  $n>2$, then
$\xi_n\leq\frac{n^2-3}{n-1}.$ 
Hence, for any $n$, we have $n\leq \xi_n<n+1$. 
Nowadays the exact values of $\xi_n$  are known for $n=2$, $n=5$, $n=9$ 
and for the infinite set of $n$ such that there exists an  Hadamard matrix of order $n+1$.
(For Hadamard matrices see, e.\,g., \cite{nu_bib_0}, \cite{nu_bib_1}.)
It is very interesting that in all these cases except $n=2$ holds the equality 
$\xi_n=n$.
Still it remains unknown whether there exists an even $n$ with such a property.
For all $n,$ if any simplex $S$ satisfies 
$S\subset Q_n\subset nS$, then the center of gravity of $S$ coincides with the center
of the cube.

In \cite{nu_bib_9} the authors introduced the notion   
{\it perfect simplex} for such a simplex $S\subset Q_n$ that $Q_n$ is inscribed 
into the simplex $\xi_n S$. 
The existence of perfect simplices earlier was known only for $n=1$ and $n=3$.
Our approach made it possible to prove that    
such simplices also exist in ${\mathbb R}^5$.

Marek Lassak 
\cite{nu_bib_2}
was the first who proved that if a simplex 
$S\subset Q_n$ 
has the maximum possible
volume, then
\begin{equation}\label{nev_ukh_eq_to1}
d_1(S)=\ldots=d_n(S)=1.
\end{equation}
Another way to get 
(\ref{nev_ukh_eq_to1})
was given in 
\cite{nu_bib_4}, where it was shown that for any maximum volume simplex 
$S\subset Q_n$ we have  $\alpha(S)=n$. In view of 
(\ref{nev_uhl_alpha_d_i_formula}) 
this is equivalent to 
(\ref{nev_ukh_eq_to1}).

By {\it $a/b$-matrix}, we mean a matrix, whose any element is equal  
to one of two numbers $a$ or $b.$ 
Values $h_n$ and $g_n$  are defined as maximum determinants of 
$0/1$ and $-1/1$-matrices of order $n$ respectively. Denote by $\nu_n$ 
the maximum volume of an $n$-dimensional simplex contained in
$Q_n$. These numbers are connected by the equalities
$g_{n+1}=2^{n}\,h_n$,
$h_n=n!\,\nu_n$, see 
\cite[Theorem 2.1]{nu_bib_1}.

If maximum $0/1$-determinant of order $n$ is known, it is possible to 
construct a maximum volume simplex in $Q_n$.
Let us enlarge the row set of such a determinant by
the row $(0,\ldots,0).$ Then simplex $S$ with these vertices 
is contained in $Q_n$ and has maximum possible volume.
Actually, consider for $S$ the node matrix ${\bf A}$ of order $n+1$, see Section 1. 
Then 
$$\vo(S)=\frac{|\det({\bf A})|}{n!}=\frac{h_n}{n!}=\nu_n.$$
Nonzero vertices of an $n$-dimensional simplex  in $Q_n$
with maximum volume one can also obtain in the same way using 
the columns of a maximum $0/1$-determinant of order $n$. 

On the other hand, if any $n$-dimensional simplex $S\subset Q_n$ 
with the zero-vertex
has maximum possible volume among all simplices in $Q_n$, 
then the coordinates of its  non-zero vertices, 
being written in rows or in columns, form maximum $0/1$-determinants of order $n.$

\section{Main Results}
\label{nev_ukh_main_results}

{\bf Theorem 1.}
{\it  Suppose {\bf M} is an arbitrary nondegenerate $0/1$-matrix of order $n$.
Let $\bf A$ be the matrix of order $n+1$ which appears
from ${\bf M}$
after adding the $(n+1)$th row $(0,0,\ldots,0,1)$ and the $(n+1)$th column
consisting of $1$'s. Denote ${\bf A}^{-1}=(l_{ij}).$ Then the following
propositions  hold true.

\smallskip
$1.$ For all $i=1,\ldots,n$,
\begin{equation}\label{nev_ukh_sum_abs_l_ij_geq2} 
\sum_{j=1}^{n+1} |l_{ij}|\geq 2. 
\end{equation}

\smallskip
$2.$ If $|\det({\bf M})|=h_n$, then 
for all $i=1,\ldots,n$ 
\begin{equation}\label{nev_ukh_sum_abs_l_ij_eq2}
 \sum_{j=1}^{n+1} |l_{ij}|= 2. 
\end{equation}

\smallskip
$3.$ If for some $i$ we have  the strong inequality
\begin{equation}\label{nev_ukh_sum_abs_l_ij_greater2}
 \sum_{j=1}^{n+1} |l_{ij}|> 2,
\end{equation}
then  $|\det({\bf M})|<h_n$. 
}

\smallskip
{\it Proof.} Consider the simplex $S$ with the zero-vertex and the rest
vertices coinciding in coordinate
form with the rows of ${\bf M}$. Obviously, 
${\bf A}$
is the node matrix of $S$. Since $\det({\bf A})=\det({\bf M})\ne 0,$
simplex $S$ is nondegenerate. Since ${\bf M}$ is a $0/1$-matrix, we have
$S\subset Q_n$. Therefore, $d_i(S)\leq 1$. 
Combined with the formula 
for the axial diameters, see (\ref{nev_uhl_d_i_formula}),
this gives 
(\ref{nev_ukh_sum_abs_l_ij_geq2}).

Now assume that $|\det({\bf M})|=h_n$. Then we have
$$\vo(S)=\frac{|\det({\bf A})|}{n!}=\frac{|\det({\bf M})|}{n!}=\frac{h_n}{n!}=\nu_n.$$
Thus, $S$ is a simplex with maximum possible volume in $Q_n.$ As it was noted above,
for such a simplex all the axial diameters $d_i(S)$ are equal to 1,
see (\ref{nev_ukh_eq_to1}).
According to  (\ref{nev_uhl_d_i_formula})
we obtain the conditions 
(\ref{nev_ukh_sum_abs_l_ij_eq2}).

Finally, if there exists $i$ such that we have  the strong inequality
(\ref{nev_ukh_sum_abs_l_ij_greater2}), then the corresponding axial diameter
satisfies $d_i(S)<1.$ Therefore, in this case the volume of $S$ is not maximum
and the absolute value of the determinant  of $0/1$-matrix ${\bf M}$ 
is strictly smaller than $h_n$. 
The theorem is proved. \hfill$\Box$

\smallskip
The inverse proposition to the second part of Theorem 1 is
not true. As a simple example, consider the case when ${\bf M}$ is the identity
matrix of order $n$. Clearly, $|\det({\bf M}|=h_n$ only for $n=1$. But since 
$$
{\bf A} =
\left( \begin{array}{ccccc}
1&0&\ldots&0&1\\
0&1&\ldots&0&1\\
\vdots&\vdots&\vdots&\vdots&\vdots\\
0&0&\ldots&1&1\\
0&0&\ldots&0&1
\end{array}
\right), \quad
{\bf A}^{-1} =
\left( \begin{array}{ccccc}
1&0&\ldots&0&-1\\
0&1&\ldots&0&-1\\
\vdots&\vdots&\vdots&\vdots&\vdots\\
0&0&\ldots&1&-1\\
0&0&\ldots&0&1
\end{array}
\right),$$
the equalities (\ref{nev_ukh_sum_abs_l_ij_eq2}) are fulfilled for all $n$ and $i=1,\ldots,n$.
The corresponding simplex $S$ is "the corner simplex" \,in $Q_n.$ Nonzero vertices  
$x^{(1)},$ 
$\ldots$ 
$x^{(n)}$ of $S$ coincide with the standard basis of ${\mathbb R}^n$ and
the last vertex $x^{(n+1)}$ is $(0,\ldots,0)$. 
Actually, though all the axial diameters $d_i(S)$ are equal to 1, the volume of $S$
becomes  maximum only in trivial case $n=1.$ 

\smallskip
Now let us note that $0/1$-matrices having maximum determinant can be obtained from maximum 
$-1/1$-determinant matrices.
Let ${\bf U}$ be an nondegenerate $-1/1$-matrix of order $n+1$ and      
${\bf T}$ be a $0/1$-matrix of order $n.$ Suppose these matrices
are connected by the following procedure considered in 
\cite[Section~2]{nu_bib_1}.  

\begin{enumerate}
\item[1.] Each column of ${\bf U}$ which starts with $-1$ 
we multiply by $-1.$
\item[2.] Each row of the new matrix which starts
with $-1$ we also multiply by
$-1.$ Denote by $\bf V$ 
the matrix of order $n+1$ coming out after steps 1 and 2.
Both the first column and the first row of this $-1/1$-matrix 
now completely consist of 1's.
\item[3.] 
Denote by 
${\bf W}$ 
the submatrix of ${\bf V}$ which
stands in rows and columns numbering $2,$ $\ldots,$ $n+1.$ 
Let us change the elements of ${\bf W}$ 
replacing $1$ by $0$ and  $-1$ by $1.$ 
By definition, ${\bf T}$ is equal to
the resulting matrix of order $n.$ 
\end{enumerate}

\smallskip
{\bf Theorem 2.}
{\it
We have the equality
\begin{equation}\label{nevskii_theor2}
|\det({\bf T})|=\frac {|\det({\bf U})|}{2^n}. 
\end{equation}
If $|\det({\bf U})|=g_{n+1},$ then    
$|\det({\bf T})|=h_n.$  
Conversely, the condition  $|\det({\bf T})|\ne h_n$ implies 
$|\det({\bf U})|\ne g_{n+1}.$  
}  

\smallskip
{\it Proof.} 
We utilize the arguments from 
\cite{nu_bib_1}.
Define simplices $D_1,$ $D_2$, and $D_3$ as follows.    
Let $D_1$ be the $(n+1)$-dimensional simplex with zero-vertex and the rest vertices
coinciding with the rows of ${\bf V}.$ 
Assume 
$D_2$ is the $n$-dimensional simplex with one vertex
$(1,\ldots, 1)$ and the rest vertices corresponding to the rows of ${\bf W}.$
Finally, let $D_3$ be the $n$-dimensional simplex whose one vertex
is $0$ and the rest are given by the rows of ${\bf T}.$ 
Clearly,
\begin{equation}\label{th2_9_2}\vo(D_1)=\frac{|\det({\bf V})|}{(n+1)!}=   
\frac{|\det({\bf U})|}{(n+1)!}. 
\end{equation}
Non-zero vertices of $D_1$ belong to the facet $x_1 = 1$ of the cube $ [-1,1]^n,$ 
therefore,
the height of $D_1$ dropped from the zero-vertex, is equal to $1.$ Further,
simplex $D_2$
is congruent to the facet of $D_1$ lying on the indicated cube facet.
It means that
\begin{equation}\label{th2_9_3}
\vo(D_1)=\frac{\vo(D_2)}{n+1}. 
\end{equation}
Changes of numbers noted in step 3 are the result of the affine
transformation of $[-1,1]^n$ into $[0,1]^n,$ for which the vertex
$(1,\ldots,1)$ of the first cube goes to the vertex $(0, \ldots, 0)$ of the second cube.
Under this transformation, measures of sets have to be multiplied by $2^{-n}$.
Since $D_3$ is the image of $D_2,$ then
\begin{equation}\label{th2_9_4}
\vo(D_3)=\frac{\vo(D_2)}{2^n}. 
\end{equation}
Moreover,
\begin{equation}\label{th2_9_5}
\vo(D_3)=\frac{|\det({\bf T})|}{n!}. 
\end{equation}
Applying
(\ref{th2_9_5}),
(\ref{th2_9_4}), and
(\ref{th2_9_3}),
we obtain
$$|\det({\bf T})|=n!\,\vo(D_3)=\frac{n!\,\vo(D_2)}{2^n}=
\frac{(n+1)!\,\vo(D_1)}{2^n}=\frac{|\det({\bf X})|}{2^n}.$$
Equality (\ref{nevskii_theor2}) is proved.
Now assume that
$|\det({\bf U})|=g_{n+1}.$ Then we have   
$|\det({\bf T})|=2^{-n}g_{n+1}=h_n.$
Obviously, the condition  $|\det({\bf T})|\ne h_n$ gives 
$|\det({\bf X})|\ne g_{n+1}.$  
This concludes the proof.
\hfill$\Box$

\smallskip
Steps 2 and 3 of the above procedure can be inversed.
Starting with a $0/1$-matrix ${\bf T}$ 
of order $n$ the inverse procedure gives  $-1/1$-matrix
${\bf V}$ 
of order $n + 1.$  
Both the first row and the first column of  ${\bf V}$ 
consist of 1's and
$|\det({\bf V})|$ $=$ $2^n|\det({\bf T})|.$ If 
$|\det({\bf T})|$ $=$ $h_n,$ then $|\det({\bf V})|$ $=$ $g_{n+1}.$ 

Having information on maximum $-1/1$-determinant of the concrete order $n$
and acting in the way described above, one can find the vertices of simplex 
$S\subset Q_n$ with maximum
volume. Values 
$\xi(S)$ 
allows to obtain rather precise estimates $\xi_n \leq \xi (S)$.
This approach was used
by the authors in \cite{nu_bib_8} for dimensions $n\leq 118$.

\begin{figure}[h]
\begin{center}
\includegraphics[width=\textwidth]{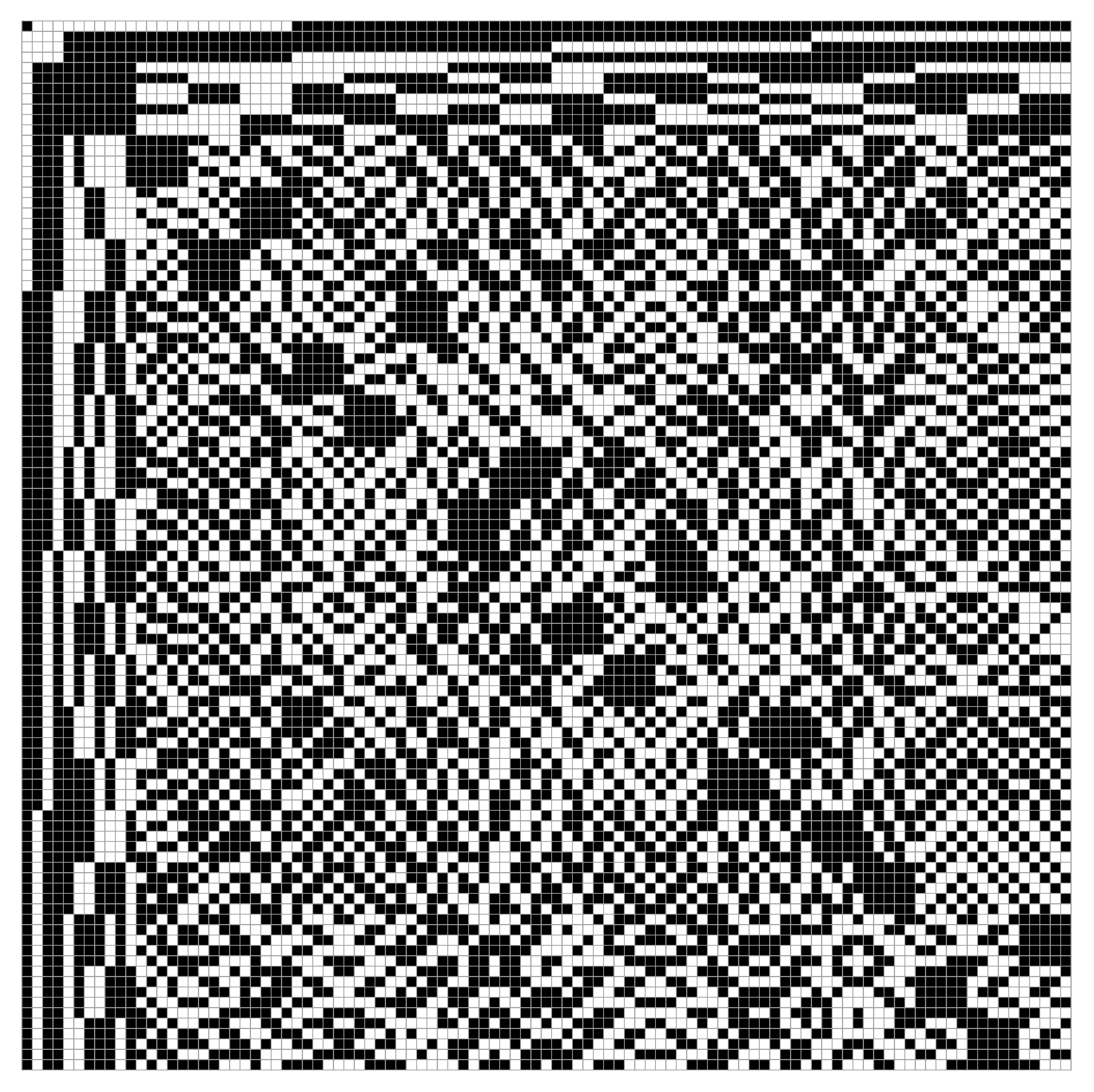}
\end{center}
\caption{$-1/1$-matrix matrix ${\bf U}$ of order $101$. White squares denote $-1$'s, black squares denote~$1$'s}
\label{fig:nev_uhl_matrix_plus_minus_ones}
\end{figure}

\section{An Example for $n=101$}
\label{nev_ukh_an_example}

Utilizing results from the previous section we will show that biggest known (by 2003) determinant of order $101$ is not maximum. For the sake of compactness, in this section we present matrices as pictures.

The collection of maximum known determinants of orders $1,\dots,119$ can be found on the site
\cite{nu_bib_12}.
We will consider the matrix of order 101 from this site. 
According the information from the page \newline
http://www.indiana.edu/\char`\~maxdet/d101.html
 this matrix was constructed by William Orrick and Bruce Solomon in~2003.  
 Note that Orrick and Solomon don't claim that this is $-1/1$-matrix 
 of largest possible determinant. 
 They only state that its determinant ``surpasses the previous record''.
\begin{figure}[h]
\begin{center}
\includegraphics[width=\textwidth]{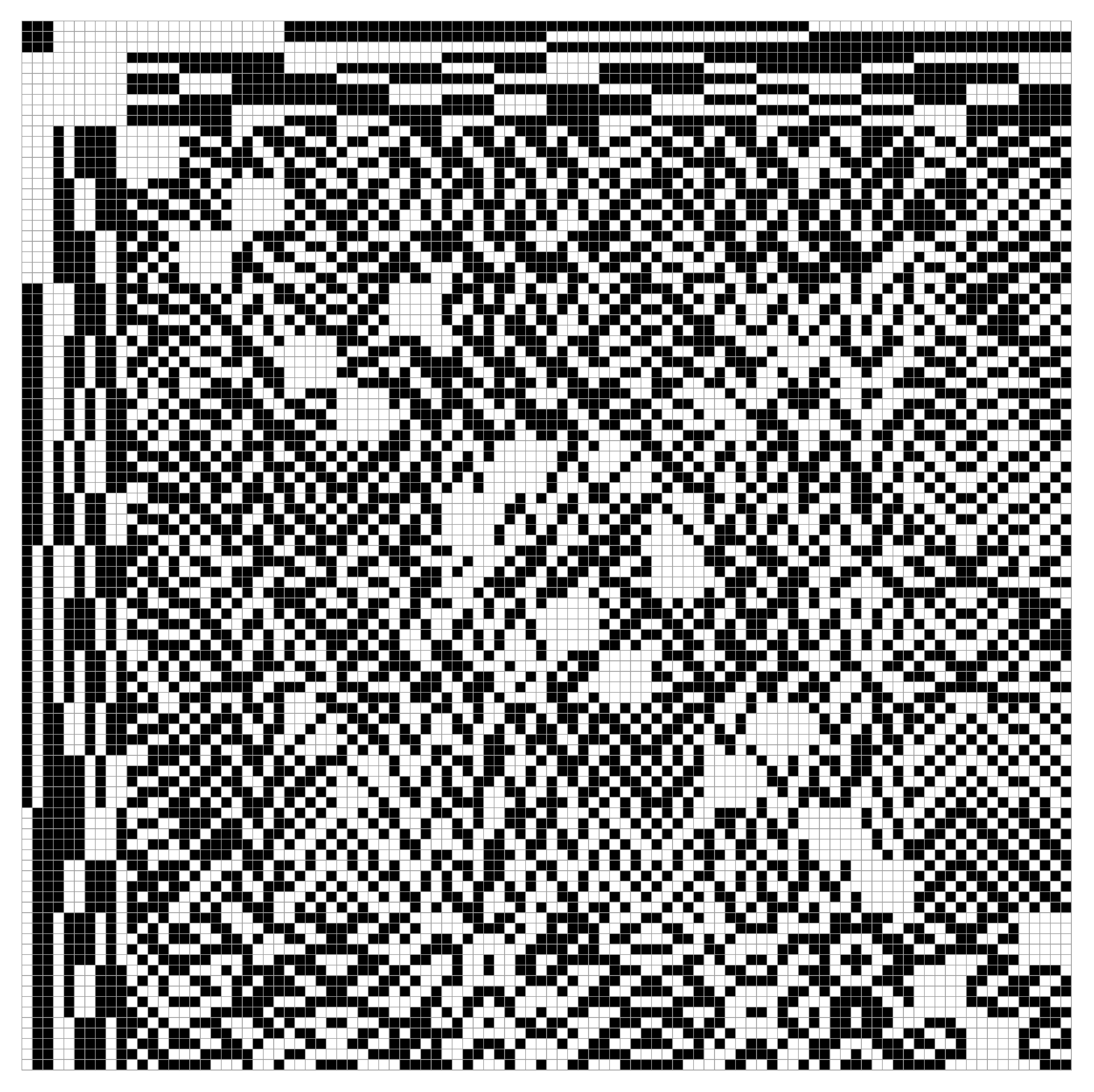}
\end{center}
\caption{$0/1$-matrix ${\bf T}$ of order $100$. White squares denote $0$'s, black squares denote $1$'s}
\label{fig:nev_uhl_matrix_zero_ones}
\end{figure}
The matrix is shown on Figure~\ref{fig:nev_uhl_matrix_plus_minus_ones},  
where white squares represent~$-1$'s while black squares denote~$1$'s. 
Using the notations of Section~\ref{nev_ukh_main_results} we will call this matrix~${\bf U}$. 
In this case $n=100$. 
Following the procedure described above, let us transform ${\bf U}$ to $0/1$-matrix ${\bf T}$ of order~$100$. 
The resulting matrix is shown on Figure~\ref{fig:nev_uhl_matrix_zero_ones}. In this case
white squares represent~$0$'s, black squares denote~$1$'s.
We will apply Theorem~1 for ${\bf M}={\bf T}$. Let us build the matrices 
${\bf A}$ and ${\bf A}^{-1}$ and calculate the sums in the left part of (\ref{nev_ukh_sum_abs_l_ij_geq2}). This gives
$$
\sum_{j=1}^{101} |l_{ij}| = 
\left\{
 \arraycolsep=3pt\def\arraystretch{1.2} 
\begin{array}{cl}
 \frac{1438}{711}=2.0225\ldots, & i=1,2,3, \\
 \frac{490}{237}=2.0675\ldots, & i=4, \\
 2, & i=5,\ldots,100. \\
\end{array}
\right. 
$$
As we see, for $i=1,2,3,4$ strict inequality \ref{nev_ukh_sum_abs_l_ij_greater2} holds.
 It follows from Theorem~1 that
$|\det({\bf T})|<h_{100}$, i.\,e., determinant of  ${\bf T}$ is not maximum $0/1$-determinant of order $100$. Consequently, by Theorem~2, determinant of matrix ${\bf U}$ is not maximum $-1/1$-determinant of order $101$.

If we consider the matrix ${\bf A}$ as the node matrix of the corresponding simplex $S$, we will find that some axial diameters of $S$ are lesser than $1$. Indeed, formula~(\ref{nev_uhl_d_i_formula}) gives
$$
d_i(S) = 
\left\{
 \arraycolsep=3pt\def\arraystretch{1.2} 
\begin{array}{cl}
\frac{711}{719}=0.9888\ldots, & i=1,2,3, \\
\frac{237}{245}=0.9673\ldots, & i=4, \\
 1, & i=5,\ldots,100. \\
\end{array}
\right. 
$$
This means that $S$ is not maximum volume simplex in $Q_{100}$.

All the calculations in this section were performed with the use of Wolfram Mathematica in symbolic mode. Accordingly, we omited most of intermediate results. Corresponding data and programs are available at\newline
 \url{http://dx.doi.org/10.17632/sm3x4xrb42.1}

\end{document}